\documentclass{article}  

\usepackage{amssymb}

\title{Decay of correlations for non H\"{o}lderian dynamics. A
  coupling approach%
\thanks{Work done within the
Projeto Tem\'atico ``Fen\^omenos Cr\'{\i}ticos em Processos Evolutivos e
Sistemas em Equil\'{\i}brio'', supported by FAPESP (grant 95/0790-1),
and is part of the activities of the N\'ucleo de Excel\^encia 
``Fen\^omenos Cr\'{\i}ticos em Probabilidade e Processos
Estoc\'asticos'' (grant 41.96.0923.00)}
}
\author{Xavier Bressaud%
\thanks{Work partially supported by FAPESP (grant 96/04860-7).}
\and
Roberto Fern\'andez%
\thanks{Researcher of the National Research Council 
(CONICET), Argentina.}\ %
\thanks{Work partially supported by CNPq (grant
301625/95-6).}% and FAPESP (grant 95/8896-3).}
\and Antonio Galves%
\thanks{Work partially supported by CNPq (grant 301301/79).}
}

\date{{\large \bf Draft} \\[10pt] 
\today}

\makeatletter
 
\@addtoreset{equation}{section}

\makeatother

\newcommand{\bit}{\begin{itemize}}
\newcommand{\eit}{\end{itemize}}
\newcommand{\bid}{\begin{description}}
\newcommand{\eid}{\end{description}}
\newcommand{\bqy}{\begin{eqnarray}}
\newcommand{\eqy}{\end{eqnarray}}
\newcommand{\beq}{\begin{equation}}
\newcommand{\eeq}{\end{equation}}
\newcommand{\pro}{{\bf P}}
\newcommand{\esp}{{\bf E}}

\newcommand{\trans}{P}

\newcommand{\transt}{{\widetilde{P}}}
\newcommand{\transb}{{\overline{P}}}

\newcommand{\atom}[3]{#1_{#2,#3}}
\newtheorem {lem}{Lemma}
\newtheorem {remark}{Remark}
\newenvironment{rem}{\begin{remark}\rm}{\end{remark}}
\newtheorem {cor}{Corollary}
\newtheorem {theorem}{Theorem}
\newtheorem {prop}{Proposition}

\newcommand{\real}{{\Bbb R}}
\newcommand{\nat}{{\Bbb N}}
\newcommand{\rel}{{\Bbb Z}}
\def\tiende#1{\mathrel{\mathop{\longrightarrow}\limits_{#1}}}
\def\var{{\rm var}}
\def\reff#1{(\ref{#1})}
\def\sqr{\vcenter{
         \hrule height.1mm
         \hbox{\vrule width.1mm height2.2mm\kern2.18mm\vrule width.1mm}
         \hrule height.1mm}}                  % This is a slimmer sqr.
\def\square{\ifmmode\sqr\else{$\sqr$}\fi}

\begin{document}          
\maketitle

\begin{abstract} We present an upper bound on the mixing rate
  of the equilibrium state of a dynamical systems defined by the
  one-sided shift and a non H\"{o}lder potential of summable
  variations.  The bound follows from an estimation of the relaxation
  speed of chains with complete connections with summable decay, which
  is obtained via a explicit coupling between pairs of chains with
  different histories.
\end{abstract}

\section{Introduction}
Let $\mu_\phi$ be the equilibrium state associated to the continuous
function $\phi$.  In this paper we obtain 
upper bounds for the speed of convergence of the limit
\beq
\label{mixing}
\int_X f \circ T^n\, g \,d\mu_\phi \tiende{n\to \infty} \int_X f \,
d\mu_\phi \int_X g \, d\mu_\phi 
\eeq 
for $\phi$ with summable
variations and $T$ the one-sided shift.  We show that this speed is
(at least) summable, polynomial or exponential according to the decay
rate of the variations of $\phi$.  The bounds apply for $f \in
L^1(\mu_\phi)$ and $g$ with variations decreasing proportionally to
those of $\phi$.
%To obtain these results, we write the left-hand side of
% (\ref{mixing}) in terms of a chain with complete connections
% [Doeblin and Fortet (1937), Lalley (1986)] and bound the speed
% of relaxation of the chain by coupling trajectories with different
% histories. 

%A dynamical system $(X,T,\mu)$ is strongly mixing if the convergence 
%occurs for all functions $f$, $g$ in a dense subset of $L^2(\mu)$. The speed
%  of this convergence ---called speed of decay of correlations, or
%  mixing rate--- is an important element in the description of the
%  dynamical system $(X,T,\mu)$. In general, it depends on the
%  regularity of the observables $f$ and $g$. 

\paragraph{}
%One-sided shifts yield simple models useful to understand a wide
%class of dynamical systems. 
%To a continuous function --- called here interaction ---  $\phi$, 
%one can associate (at least) an
%equilibrium state $\mu_\phi$. The mixing properties of the one-sided shift
%with respect to this measure depend strongly on the regularity properties of
%$\phi$.  It was proved by Bowen (1975) that such an equlibrium state is unique
%if $\phi$ has summable decay. 
Previous approaches to the study of the mixing properties of the
one-sided shift rely on the use of 
the transfer operator  $L_\phi$, % associated to the interaction $\phi$ 
defined by the duality, 
\beq
\label{duality}
\int_X f \circ T^n \,g\, d\mu_\phi = \int_X f \,L_\phi^n g\,
d\mu_\phi\;.  \eeq 
If $\phi$ is H\"{o}lder, this operator, acting on the subspace of
H\"older observables, has a spectral gap and the limit (\ref{mixing})
is attained at exponential speed (Bowen, 1975).  When $\phi$ is not
H\"{o}lder, the spectral gap of the transfer operator may vanish and
the spectral study becomes rather complicated.  To estimate the mixing
rate, Kondah, Maume and Schmitt (1996) proved first that the operator
is contracting in the Birkhoff projective metric, while Pollicott
(1997), following Liverani (1995), considered the transfer operator
composed with conditional expectations.  In contrast, our approach is
based on a probabilistic interpretation of the duality~(\ref{duality})
in terms of expectations, conditioned with respect to the past, of a
chain with complete connections The convergence (\ref{mixing}) is
therefore related to the relaxation properties of this chain.  In this
paper, such relaxation is studied via a coupling method.

\paragraph{} 
Coupling ideas were first introduced by Doeblin in his 1938 work on
the convergence to equilibrium of Markov chains.  He let two
independent trajectories evolve simultaneously, one
starting from the stationary measure and the other from an arbitrary
distribution.  The convergence follows from the fact that both
realizations meet at a finite time.  Instead of letting the
trajectories evolve independently, one can couple them from the
beginning, reducing the ``meeting time'' and, hence, obtaining a
better rate of convergence (leading to the so-called Dobrushin's
ergodic coefficient).  Doeblin published his results in a hardly known
paper in the Revue Math\'ematique de l'Union Interbalkanique. (For a
description of Doeblin's contributions to probability theory we refer
the reader to Lindvall 1991). His ideas were taken up and exploited
only much later in papers by Athreya, Ney, Harris, Spitzer and
Toom among others.  
The sharpness of the convergence rates provided by different types of
Markovian couplings has been recently discussed by Burdzy and Kendall
(1998). 

\paragraph{}
In the context of dynamical systems, the recent papers by Coelho and
Collet (1995) and Young (1997) consider the time two independent
systems take to become close.  This is reminiscent of the original
coupling by Doeblin.  In Bressaud, Fern\'andez, Galves (1997), the
coupling approach was generalized to treat chains with complete
connections.  These processes, introduced by Doeblin and Fortet (1937)
(see also Lalley, 1986) appear in a natural way in the context of
dynamical systems.  They are characterized by having transition
probabilities that depend on the whole past, albeit in a continuous
manner.  Due to this fact, the coupling can not ensure that two
different trajectories will remain equal after their first meeting
time.  But the coupling used in the present paper, and in our
preceeding one, has the property that if the trajectories meet they
have a large probability of remaining equal, and this probability
increases with the number of consecutive agreements.  In the summable
case, the coupling is such that with probability one the trajectories
disagree only a finite number of times.  In fact, our approach can
also be applied under an assumption weaker than summability
[\reff{weak.1} below], leading to trajectories that differ infinitely
often but with a probability of disagreement that goes to zero.  The
method leads, in particular, to a criterium of uniqueness for
$g$-measures proven by Berbee (1987).  The mean time between succesive
disagreements provides a bound on the speed of relaxation of the chain
and hence, through our probabilistic interpretation of
(\ref{duality}), of the mixing rate.
%This coupling yields an estimate of the relaxation speed of the chains with
%complete connections involved in~(\ref{dualityproba}) and, hence, a
%bound on the speed of decay of correlations for the corresponding equlibrium
%state for different types of observables.  

%We point out that Ferrari, Maas and Martinez (1998), in their studies
%of regenerative representations of chains with complete connections,
%introduced couplings that share many properties with ours.  The
%relations between their coupling and ours deserves further
%investigation.  
%We believe that the coupling approach has a promisory potential and
%that our method in particular could be useful in other
%settings as well.  On the one hand, the multidimensional version of
%our coupling leads to a criterium of uniqueness of Gibbs states, which
%extends to long-range interactions the criterium based on disagreement
%percolation developed by van den Berg (1992) and van den Berg and Maes
%(1992).  On the other hand, we expect our methods to be also
%applicable when the spectral gap of the transfer operator vanishes due
%to lack of hyperbolicity (eg.\ for the so-called intermitent maps).
%Previous studies of this problem have resorted to Birkhoff projective
%metric [Saussol, Vaienti, Liverani (1997)] and to properties of the
%associated zeta-function [Isola (1997)].  Coupling ideas have also
%been recently applied by Young (1997).

\paragraph{}
The paper is organized as follows. The main results and definitions
relevant to dynamical systems are stated in Section~\ref{definitions}.
The relation between chains with complete connections and the transfer
operator is spelled out in Section~\ref{transferoperator}.  In
Section~\ref{relaxation}, we state and prove the central result on
relaxation speeds of chains with complete connections.
Theorem~\ref{normalizedth} on mixing rates for normalized functions
is proven in Section~\ref{normalized}, while Theorem~\ref{general} on
rates for the general case is proven in Section~\ref{proof}.  The
upper bounds on the decay of correlations depend crucially on
estimations of the probability of return to the origin of an auxiliary
Markov chain, which are presented in Section~\ref{markovchain}.

%The mixing rates therein are determined by the return times of a very
%simple Markov chain [defined in (\ref{rr.1})--(\ref{rr.5}) below].
%Nevertheless, there seems to be no result on these times in the
%literature.  We present some estimations in Section~\ref{markovchain}.

\section{Definitions and statement of the results}
\label{definitions}

\paragraph{}
Let $A$ be a finite set henceforth called \emph{alphabet}. Let us
denote 
\begin{equation}
  \label{eq:40}
  \underline{A} \;=\; \bigl\{x = (x_j)_{j\leq-1}\;,\; x\in A\bigr\}
\end{equation}
the set of sequences of elements of the alphabet indexed by the
strictly negative integers.  Each sequence $x\in\underline A$ will be
called a {\em history}.  Given two histories $x$ and
$y$, the notation $x \stackrel{m}{=}y$ indicates that $x_j = y_j$ for
all $-m \leq j \leq -1$.

As usual, we endow the set $\underline A$ 
with the product topology and the $\sigma$-algebra generated by the
cylinder sets.  We denote
by ${\cal C}^0(\underline{A}, \real)$ the space of real-valued
continuous functions on $\underline{A}$. 

\paragraph{}
We consider the one-sided shift $T$ on $\underline{A}$, 
$$ \begin{array}{llcl}
T:& \underline{A} & \longrightarrow & \underline{A} \\
& x& \longmapsto & T(x) = (x_{i-1})_{i\leq -1}. 
\end{array}
$$
Given an element $a$ in $A$ and an element $x$ in $\underline{A}$,
we shall denote by $xa$ the
element $z$ in $\underline{A}$ such that $z_{-1} = a$ and
$T(z) = x$. 

\subparagraph{}
Given a function $\phi$ on $\underline{A}$, 
$ \phi:  \underline{A}  \to  \real$, 
we define its %{\em regularity rate} to be the
sequence of \emph{variations} $(\var_{m}(\phi))_{m \in \nat}$, 
\beq
\label{defgamma}
\var_{m}(\phi) = \sup_{x \stackrel{m}{=}y}|\phi(x) - \phi(y)| \;.
\eeq
%
%Such a function $\phi$ can be interpreted as an energy density; we
%shall call it an {\em interaction} 
We shall say that it has {\em summable variations} if,
\begin{equation}
  \label{rr.10}
\sum_{m\geq 1} \var_m(\phi) < + \infty\;,   
\end{equation}
and that it is {\em normalized} if it satisfies, 
\begin{equation}
  \label{rr.15}
  \forall x \in \underline{A}\;, \;\;\; \sum_{a \in A} e^{\phi(xa)} =
1\;.
\end{equation}
We say that a shift-invariant measure $\mu$ on $\underline{A}$
is {\em compatible with the normalized function $\phi$} if
and only if, for $\mu_\phi$-almost-all $x$ in $\underline{A}$,
\begin{equation}
  \label{eq:50}
\esp_{\mu_\phi}\bigl( 1_{\{x_{-1}=a\}}|{\cal F}_{\le -2}\bigr)(x) 
\; =\; e^{\phi(T(x)a)}\;, 
\end{equation}
where the left-hand side is the usual conditional expectation of the
the indicator function of the event $\{x_{-1}=a\}$ with respect to the
$\sigma$-algebra of the past up to time $-2$.  

An equivalent way of expressing this is by saying that $\mu_\phi$ is a
$g$-measure for $g = e^\phi$. If $\phi$ has summable variations, and
even under a slightly weaker conditions, then such a measure is unique
and will be denoted $\mu_\phi$.  The measure $\mu_\phi$ can also be
characterized via a variational principle, in which context it is
called \emph{equilibrium state} for $\phi$.  For details see
Ledrappier (1974), Walters (1975), Quas (1996) and Berbee (1987).

\paragraph{}
For a non-constant $\phi$, we consider the seminorm
\begin{equation}
  ||g||_{\phi} \;=\; \sup_{k \geq 0}
\frac{ \var_k(g)}{\var_k(\phi)}
\end{equation}
and the subspace of ${\cal C}^0(\underline{A})$ defined by,
\begin{equation}
  \label{rr.20}
  V_{\phi} \;=\; \biggl\{ g \in {\cal C}^0(\underline{A}, \real)\,, \;
 ||g||_\phi < + \infty \biggr\} \;.
\end{equation}

\paragraph{}
Given a real-valued sequence $(\gamma_{n})_{n \in \nat}$, let
$(S^{(\gamma)}_{n})_{n \in \nat}$ be the Markov chain taking values
in the set $\nat$ of natural numbers starting from the origin
\begin{equation}
\pro(S^{(\gamma)}_0 =0)\;=\;1
\label{rr.5.-1}
\end{equation}
whose transition probabilities are defined by
\begin{equation}
  \label{rr.5}
\begin{array}{rcl}
p_{i,i+1} &=& 1- \gamma_i \\
p_{i,0} &=& \gamma_i \;,
\end{array}
\end{equation}
for all $i\in\nat$.
For any $n\ge 1$ we define %the sequence $(\gamma^*_{n})_{n \in \nat}$ by 
\begin{equation}
\gamma^*_n = \pro(S^{(\gamma)}_n = 0)\;.
\label{rr.1}
\end{equation}

\paragraph{}
We now state our first result. 
\begin{theorem}
\label{normalizedth}
Let $\phi: \underline{A} \to \real$ be a normalized function with summable
variations and set 
\begin{equation}
  \label{rr.30}
  \gamma_n \;=\; 1 - e^{-\var_n(\phi)} \;.
\end{equation}
Then, 
\begin{eqnarray}
  \label{rr.35}
    \left| \int f \circ T^n \,g\, d\mu_\phi - \int f\,
  d\mu_\phi \int g \, d\mu_\phi \right| 
&\leq& \,  ||f||_1\,   ||g||_{\phi} \, \sum_{k=0}^n \var_k(\phi)\,
\gamma^*_{n-k} \label{eq:x}\\[5pt]
&\leq& C \,  ||f||_1\,   ||g||_{\phi} \, \gamma^*_{n}\;, 
\label{eq:xx}
\end{eqnarray}
for all $f \in L^1(\mu_\phi)$ and $g \in V_{\phi}$, for a computable
constant $C$.
\end{theorem}
This theorem is proven in Section \ref{normalized}, using the results
obtained in Section \ref{relaxation} on the relaxation speed of chains
with complete connections.

\paragraph{}
For each non-normalized function $\phi$ with summable variations there
exists a unique positive function $\rho$ such that the function
\begin{equation}
  \label{rr.25}
  \psi = \phi + \log{\rho} - \log{\rho \circ T}
\end{equation}
is normalized (Walters, 1975). We call $\psi$ the {\em normalization}
of $\phi$.  The construction of compatible measures given in
\reff{eq:50} looses its meaning for non-normalized $\phi$.  It is
necessary to resort to an alternative characterization in terms of a
variational principle (see eg.\ Bowen 1975) leading to equilibrium
states.  In Walters (1975) it is proven that: 
\begin{itemize}
\item[(a)] $\phi$ with
summable variations admits a unique equilibrium state, that we denote
also $\mu_\phi$; 
\item[(b)] the corresponding normalized $\psi$, given by
\reff{rr.25}, admits a unique compatible measure $\mu_\psi$ (even when
the variations of $\psi$ may not be summable), and 
\item[(c)]  $\mu_\phi=\mu_\psi$.
\end{itemize}

\paragraph{}
Our second theorem generalizes Theorem~\ref{normalizedth} to 
non-normalized functions. 
\begin{theorem}
\label{general}
Let $\phi: \underline{A} \to \real$ be a function with summable
variations and let $\psi$ be its normalization.  Let $(n_m)_{m \in \nat}$
be an increasing subadditive sequence such that the subsequence of the
rests, $\bigl(\sum_{k \geq n_{m}} \var_k(\phi)\bigr)_{m \geq 0}$, is
summable, and
\begin{equation}
  \label{rr.40}
  \overline{\gamma}_m = 1 - e^{- 3 \sum_{k \geq n_{m}}  \var_k(\phi)}\;; 
\end{equation}
then, 
\begin{eqnarray}
  \label{rr.45}
   \left| \int f \circ T^n \,g\, d\mu_\phi - \int f\,
  d\mu_\phi \int g\,  d\mu_\phi \right| &\le&
||f||_1\,   ||g||_{\phi} \, \sum_{k=0}^n \var_{n_k} (\phi)\,
\overline\gamma^*_{n-k} \\[5pt]
&\leq & C \,  ||f||_1\,   ||g||_{\phi} \, \overline{\gamma}^*_{n} \;, 
\end{eqnarray}
for all  $f \in L^1(\mu_\phi)$ and $g \in V_{\phi}$, for a computable
constant $C$.  Here
$\overline\gamma^*$ is defined as in \reff{rr.1} but using the
sequence $(\overline\gamma^*_{n})_{n \in \nat}$.
\end{theorem}

\paragraph{}
The estimation of the large-$n$ behavior of the
sequence $(\gamma^*_{n})_{n \in \nat}$ given the behavior of the
original $(\gamma_{n})_{n \in \nat}$ only requires elementary
computations. For the convenience of the reader we summarize some
results in Appendix~\ref{markovchain}. 

%\paragraph{}
%It appears from Theorem~\ref{normalizedth} that 

%\begin{cor}
%If $(\var_n(\phi))_{n \in \nat}$ has polynomial decay, then so does
%$(\gamma^*_{n})_{n \in \nat}$. More precisely, in this case, there is a
%constant $C^*$ such that $\gamma^*_{n} \leq C^* \sum_{k \geq n} \var_k(\phi)$. 
%In particular, if  there is a $r>1$ and a
%constant $C$ such that $\var_n(\phi) \leq C \frac{1}{n^r}$, then there is a
%constant $C^*$ such that $\gamma^*_{n} \leq C^* \frac{1}{n^{r-1}}$. 
%If $(\var_n(\phi))_{n \in \nat}$ decay exponentially, then so does
%$(\gamma^*_{n})_{n \in \nat}$ althought at a possibly weaker rate of decay. 
%\end{prop}
%\paragraph{Comparison with previous results.}

\section{Transfer operators and chains.}
\label{transferoperator} 

\paragraph{}
Let $P$ be a family of transition probabilities on $A
\times \underline{A}$,  
\begin{equation}
\begin{array}{llcl}
\trans:& A \times \underline{A} & \longrightarrow & [0;1] \\
& (a,z)& \longmapsto & \trans(a\,|\,z)\;.  
\end{array}
\label{eq:-1}
\end{equation}
Given a history $x$, a \emph{chain with past $x$ and transitions
  $\trans$}, is the process $(Z^{x}_{n})_{n  \in{\bf Z}}$ whose conditional
probabilities satisfy
\begin{equation}
  \label{eq:10}
  \pro(Z^{x}_{n}= a \,|\,
Z^{x}_{n+j} = z_j, j\leq -1) \;=\; \trans(a \,|\,z)
\hbox{ for } n\geq 0 \;,
\end{equation}
for all $a \in A$ and all histories $z$ with $z_{j-n}
= x_j, j\leq -1$, and such that 
\begin{equation}
  \label{eq:mais}
Z^{x}_{n}=x_n \;,\;  \hbox{ for } n\leq -1\;.
\end{equation}
This chain can be interpreted as a conditioned version of the process
defined by the transition probabilities \reff{eq:-1}, given a past $x$
(for more details, see Quas 1996).

\paragraph{}
Let $\phi: \underline{A} \to \real$ be a continuous normalized
function.  The transfer operator associated to $\phi$ is the
operator $L_\phi$ acting on ${\cal C}^0(\underline{A}, \real)$ defined
by,
\begin{equation}
  \label{eq:100}
   L_\phi f (x) = \sum_{y\, :\, T(y) =
  x} e^{\phi(y)} f(y)\;.
\end{equation}
This operator is related to the conditional probability \reff{eq:50}
  in the form
  \begin{equation}
    \label{eq:105}
    \esp_{\mu_\phi}\bigl( f \,|\,{\cal F}_{\le -2}\bigr) \;=\; 
\bigl(L_\phi f\bigr) \circ T\;.
  \end{equation}
This relation shows the equivalence of \reff{duality} and
\reff{eq:100} as definitions of the operator.
In addition, if $\phi$ is normalized we can construct, for each
history $x \in \underline{A}$, the chain
$Z^{x}_{\phi}=(Z^{x}_{n})_{n \in{\bf Z}}$ with past $x$ and transition
probabilities 
\begin{equation}
  \label{eq:20}
  \trans(a\, |\, x ) = e^{\phi(xa)}\;.
\end{equation}

\paragraph{}
Iterates of the transfer operator, $L^n_\phi g (x)$, on functions 
$g \in {\cal C}^0(\underline{A})$ can be interpreted
as expectations $\esp[g((Z^{x}_{n+j})_{j \leq -1})]$ of the chain.
Indeed, 
\bqy
L^n_\phi g (x) &=& \sum_{a_1, \ldots, a_n \in A}e^{\sum_{k = 1}^n
  \phi(x a_1 \cdots a_k )} g(x a_1 \cdots a_n)
\nonumber \\
&=& \sum_{a_1, \ldots, a_n \in A} \left( \prod_{k=1}^n
  \trans(a_k\,|\,a_{k-1}
  \cdots a_1 x) \right) g(x a_1 \cdots a_n) \nonumber \\
&=& \esp[ g((Z^{x}_{n+j})_{j \leq -1})]\;. \nonumber 
\eqy

\paragraph{}
From this expression and the classical duality~(\ref{duality}) between
the composition by the shift and the transfer operator $L_\phi$ in
$L^2(\mu_\phi)$, we obtain the following expression for the decay of
correlations, 
\begin{eqnarray}
\lefteqn{ \int f \circ T^n \,g\, d\mu_\phi - \int f\, d\mu_\phi
  \int g\,   d\mu_\phi }\nonumber\\
 &= &  \int f(x)\, L_{\phi}^n g(x)\,
  d\mu_\phi(x) - \int f(x) \left( \int L_{\phi}^n g(y) \,d\mu_\phi(y)
  \right) d\mu_\phi(x)  \nonumber \\
& =&  \int f(x) \int \Bigl(
    \esp[g((Z^{x}_{n+j})_{j \leq -1})] - \esp[g((Z^{y}_{n+j})_{j \leq
      -1})]\Bigr)\, d\mu_\phi(y)\, d\mu_\phi(x) \;.
\label{interpretation}
\end{eqnarray}
This inequality shows how the speed of decay of correlations can be
bounded by the speed with which the chain loosses its memory.  We deal
with the later problem in the next section.

\section{Relaxation speed for chains with complete connections}
\label{relaxation}
\subsection{Definitions and main result}
\paragraph{}

We consider chains whose transition probabilities satisfy 
\beq
\label{p.3}
\inf_{  x , y : x \stackrel{m}{=}y } 
{\trans(a \,|\, x) \over \trans(a \,|\, y)}  \geq \; 1 -   \gamma_m \;,
\eeq
%\paragraph{}
for some real-valued sequence $(\gamma_m)_{m \in \nat}$,
decreasing to $0$ as $m$ tends to $+ \infty$.   Without loss of
generality, this decrease can be assumed to be monotonic.  To avoid
trivialities we assume $\gamma_0<1$.  In the literature, a
\emph{stationary} process satisfying \reff{p.3} is called a
\emph{chain with complete connections}.

\paragraph{}

\paragraph{}
For a set of transition probabilities satisfying \reff{p.3}, we
consider, for each $x\in\underline A$, the chain 
$(Z^{x}_n)_{n\in\rel}$ with past $x$ and transitions $\trans$ [see
\reff{eq:10}--\reff{eq:mais}].  The following proposition
  plays a central role in the proof of our results.
\begin{prop}
\label{existcoupling}
For all histories
$x,y \in \underline{A}$, there is a coupling 
$\bigl((\widetilde{U}^{x,y}_n,\widetilde{V}^{x,y}_n)\bigr)_{n\in\rel}$
of $(Z^{x}_n)_{n\in\rel}$ and $(Z^{y}_n)_{n\in\rel}$ such that the
    integer-valued  process $(T^{x,y}_n)_{n\in {\bf Z}}$  defined by 
\beq
\label{defcompteur}
T^{x,y}_n \;=\; \inf \{ m \geq 0 \, : \, \widetilde{U}^{x,y}_{n-m} \neq
\widetilde{V}^{x,y}_{n-m} \}, 
\eeq 
satisfies 
\beq
\label{compteur}
\pro(T^{x,y}_n = 0) \;\leq \;  \gamma^*_n
\eeq
for $n\ge 0$, where $\gamma^*_n$ was defined in \reff{rr.1}.
\end{prop}
The proof of this proposition is given in Section \ref{sepro}.

\paragraph{}
An immediate consequence of this proposition is the following bound on
the relaxation rate of the processes $Z^{x}$. 
\begin{cor}%[Relaxation speed]
\label{relaxccc}
For all histories $x$ and $y$, for all $a \in A$,  
\begin{equation}
\label{le.10}
 \Bigl| \pro(Z^{x}_{n}= a) - \pro(Z^{y}_{n} = a) \Bigr| \;\leq\;
\gamma^*_n\;,
\end{equation}
and, for $k\geq 1$, 
\begin{eqnarray}
\lefteqn{\hspace{-3cm}\left| \pro\Bigl((Z^{x}_{n}, \ldots, Z^{x}_{n+k})= (a_0,
  \ldots, a_k) \Bigr) - \pro\Bigl((Z^{y}_{n}, \ldots, Z^{y}_{n+k})= (a_0,
  \ldots, a_k) \Bigr) \right|} \nonumber\\
&& \qquad\qquad \le \quad\sum_{j=0}^{k}   \left(\prod_{m=1}^{j-1}(1- \gamma_{m})
  \right) \gamma^*_{n-j}\;.
\label{le.7} 
\end{eqnarray}
\end{cor}
This lemma is proved in Section~\ref{proofcor}.

\begin{rem}
\label{unicitelarge}
Whenever 
\begin{equation}
  \label{eq:pp}
 \gamma^*_n\to 0 \;,
\end{equation}
inequality \reff{le.10} implies the
existence and uniqueness of the invariant
measure compatible with a system of conditional probabilities
satisfying \reff{p.3}.  In fact, property \reff{eq:pp} holds under
the condition 
\begin{equation}
  \label{weak.1}
\sum_{m\geq 1} \prod_{k=0}^{m}(1 - \gamma_k) \;=\;  + \infty \;. ,
\end{equation}
which is weaker than summability.
In this case, the Markov chain $(S^{(\gamma)}_n)_{n \in \nat}$ is no
longer transient but it is null recurrent and the property
$\pro(S^{(\gamma)}_n =0) \to 0$ remains true.  %Our argument to
%estimate the speed, however, must be adapted.

\end{rem}

\begin{rem}
If $X = (X_n)_{n\in {\bf Z}}$ is a stationary process with
transition $\trans$ satisfying \reff{p.3}, then
Corollary~\ref{relaxccc} implies 
\begin{equation}
 \Bigl| \pro(Z^{x}_{n}= a) - \pro(X_{n} = a) \Bigr| \;\leq\;
 \gamma^*_n\;, 
\label{le.12}
\end{equation}
uniformly in the history $x$.
\end{rem}

\subsection{Maximal coupling}
Given two probability distributions $\mu = (\mu(a))_{a \in A}$ and
$\nu = (\nu(a))_{a \in A}$ we denote by $\mu \tilde{\times} \nu = (\mu
\tilde{\times} \nu (a,b))_{(a,b) \in A \times A}$ the so-called {\em
  maximal coupling} of the distributions $\mu$ and $\nu$ defined
as follows:
\begin{equation}
\left\{ 
\begin{array}{ll}
\mu \tilde{\times} \nu (a,a) = \mu(a) \wedge \nu(a) & \hbox{ if } a
= b \\[15pt]
\displaystyle \mu \tilde{\times} \nu (a,b) = \frac{(\mu(a) -\nu(a))^+
  (\nu(b) -\mu(b))^+}{\sum_{e \in A}(\mu(e) -\nu(e))^+} 
 & \hbox{ if } a
\neq b \;. 
\end{array}
\right.  
\end{equation}
For more details on maximal couplings see Appendix A.1 in
Barbour, Holst and Janson (1992).
\paragraph{}
The coupling is maximal in the sense that the distribution $\mu
\tilde{\times} \nu $ on $A \times A$ maximizes the weight 
$$\Delta(\zeta) = \sum_{a \in A} \zeta(a,a)$$
of the diagonal among the distributions $\zeta$ on $A
\times A$ satisfying simultaneously 
$$\sum_{a \in A} \zeta (a,b) = \nu(b) \;\;\; \hbox{ and }\;\;\; \sum_{b \in A}
\zeta(a,b) = \mu(a)\;.$$  
For this coupling, the weight $\Delta(\mu \tilde{\times} \nu)$ of the diagonal satisfies, 
\beq
\label{coupling}
\Delta(\mu \tilde{\times} \nu) =   \sum_{a \in A} \mu(a) \wedge
\nu(a) = 1 - \sum_{a \in A} (\mu(a) - \nu(a))^+ = 1 - \frac{1}{2} \sum_{a
  \in A} |\mu(a) - \nu(a)|. 
\eeq
Moreover, 
\beq
\label{coupling2}
\Delta(\mu \tilde{\times} \nu) = 1 - \sum_{a \in A} \mu(a)\left(1 -
\frac{\nu(a)}{\mu(a)}\right)^+ \geq 1 - \sum_{a \in A} \mu(a)\left(1 -
\inf_{a' \in A}\frac{\nu(a')}{\mu(a')}\right) = \inf_{a \in
A}\frac{\nu(a)}{\mu(a)}.  \eeq
%and 
%\beq
%\label{coupling3}
%\Delta(\mu \tilde{\times} \nu) =  \sum_{a \in A} \mu(a) \wedge
%\nu(a) \geq \inf_{a \in A} \min (\mu(a) ; \nu(a)). 
%\eeq

\subsection{Coupling of chains with different pasts}
\paragraph{}
Given a double history $(x, y)$, we consider the transition
probabilities defined by the maximal coupling
\begin{equation}
\transt((a,b) \, | \, x,y)  = \left[\trans(\cdot \,|\,
x) \tilde{\times} \trans(\cdot  \,|\, y)\right]  (a,b)\;.
\label{ant.1}
\end{equation}
%and let $m_0$ denote the first integer for which $\gamma_{n} \leq
%\gamma_{0}$. 
By \reff{p.3} we have, 
$$ \inf_{a \in A, u \stackrel{m}{=}
 v } \frac{\trans(a\,|\,u)}{\trans(a\,|\,v)} \;\geq\;
 1 - \gamma_m. $$
%  \;\;\; \hbox{ and } \;\;\; \inf_{a \in A, x \in \underline{A}} \trans(a|x)
%\geq 1 - \gamma_0$$
By (\ref{coupling2}) this implies that 
\begin{equation}
\label{diagonale}
\Delta\left(\transt(\,\cdot\,,\,\cdot \, | \, x,y) \right)\; \geq \; 1
- \gamma_m\;,
\end{equation}
whenever $ x \stackrel{m}{=} y$.

\paragraph{}
Now, we fix a double history $(x, y)$ and we define
$\bigl((\widetilde{U}^{x, y}_n, \widetilde{V}^{x,
  y}_n)\bigr)_{n\in\rel}$ to be the chain taking values in $A^2$, with
past $(x, y)$ and transition probabilities given by \reff{ant.1}. If $
x \stackrel{m}{=} y$, (\ref{diagonale}) yields
\begin{equation}
\label{alpha}
\pro( \widetilde{U}^{x,y}_0 \neq \widetilde{V}^{x,y}_0) \;\leq\;
\gamma_{m} . 
\end{equation}
\paragraph{}
We denote 
\begin{equation}
\Delta_{m,n} \;:=\;
% \{ \widetilde{U}_j =  \widetilde{V}_j \} \;=\; 
\Bigl\{ \widetilde{U}_j =  \widetilde{V}_j \, ,\, m \leq j \leq n
\Bigr\}\;. 
\label{le.15}
\end{equation}
Notice that $\Delta_{-m,-1}$ is the reunion over all the sequences $x,
 y$ with $ x \stackrel{m}{=} y$ of the events $\{ (\widetilde{U}_j
 ,\widetilde{V}_j) = (x_j,y_j) \,;\, j\leq -1 \}$.  Using the
 stationarity of the conditional probabilities, we obtain
\begin{equation}
\label{beta}
\pro( \widetilde{U}_n \neq \widetilde{V}_n \, |\, \Delta_{n-m,n-1}) 
\;\leq\; \gamma_{m} \;,
\end{equation}
for all $n\ge 0$.

\subsection{Proof of Proposition~\protect{\ref{existcoupling}}}
\label{sepro}

From this subsection on, will be working with bounds which are uniform
in $x,y$, hence we will omit, with a few exceptions, the superscript
$x,y$ in the processes $T^{x,y}_n$ (defined below), $\widetilde
U^{x,y}_n$ and $\widetilde V^{x,y}_n$.

Let us consider the integer-valued 
process $(T_n)_{n\in {\bf Z}}$  defined by: 
\begin{equation}
T_n \;=\; \inf \{ m \geq 0 \, : \, \widetilde{U}_{n-m} \neq
\widetilde{V}_{n-m} \}\;.
\label{le.19}
\end{equation}
For each time $n$, the random variable $T_n$ counts  the number of
steps backwards needed to find a difference in the coupling.
First, notice that (\ref{beta}) implies that, 
\begin{equation}
\pro(T_{n+1}=k+1 \,|\, T_{n} = k) \;\geq\; 1-\gamma_{k}
\label{le.20}
\end{equation}
and
\begin{equation}
\pro(T_{n+1}= 0 \,|\, T_{n} = k) \;\leq\; \gamma_{k}\;,
\label{le.25}
\end{equation}
all the other transition probabilities being zero.
This process $(T_n)_{n\in {\bf Z}}$ is not a Markov chain.  

\paragraph{}
We now consider the integer-valued Markov chain $(S^{(\gamma)}_n)_{n
  \geq 0}$ starting from state $0$ and with transition probabilities
given by \reff{rr.5}, that is $p_{i,i+1} = 1- \gamma_i$ and $p_{i,0} =
\gamma_i$.  Proposition~\ref{existcoupling} follows from the following
lemma, setting $k=1$.

\begin{lem}
\label{dominlem}
For each $k \in {\bf N}$, the following inequality holds:

\beq
\label{domin}
 \pro(S^{(\gamma)}_n \geq k) \;\leq\;
\pro(T_n \geq k) 
\eeq
\end{lem}
\paragraph{Proof}
We shall proceed by induction on $n$.  Since $\pro(S^{(\gamma)}_0
=0)=1$, inequalities \reff{domin} holds for $n=0$.  Assume now that
\reff{domin} holds for some integer $n$.  There is nothing to prove
for $k=0$.. For $k \geq 1$, 
\bqy
\pro(T_{n+1} \geq k)& = &
\sum_{m=k}^{+\infty} \pro(T_{n+1} = m )\nonumber \\
%& = &
%\sum_{m=k}^{+\infty} \pro(T_{n+1} = m, T_{n} = m-1 )\nonumber \\
& = &
\sum_{m=k}^{+\infty} \pro(T_{n+1} = m \,|\, T_{n} = m-1 ) 
\,\pro(T_{n} = m-1 ) \nonumber \\
&\geq  &
\sum_{m=k}^{+\infty} (1-\gamma_{m-1})\, \pro(T_{n} = m-1 ) \nonumber \\
&=  &
\sum_{m=k}^{+\infty} (1-\gamma_{m-1}) \,\Bigl(  \pro(T_{n} \geq  m-1
)-\pro(T_{n} \geq  m )\Bigr)  \nonumber \\ 
&=  & (1-\gamma_{k-1})\,  \pro(T_{n} \geq  k-1 ) + 
\sum_{m=k}^{+\infty} (\gamma_{m-1} - \gamma_{m}) 
\,\pro(T_{n} \geq  m ) \;.
\eqy
By the same computation, we see that 
\beq
\label{recurmarkovchain}
\pro(S^{(\gamma)}_{n+1} \geq k) \;= \;
 (1-\gamma_{k-1}) \,\pro(S^{(\gamma)}_{n} \geq  k-1 ) + 
\sum_{m=k}^{+\infty} (\gamma_{m-1} - \gamma_{m})
\,\pro(S^{(\gamma)}_{n} \geq  m ) \;.
\eeq
Hence,  using  the recurrence assumption and the fact that
 $(\gamma_{n})_{n\geq 0}$ is decreasing we conclude that
$$\pro(T_{n+1} \geq k) \;\geq\; \pro(S^{(\gamma)}_{n+1} \geq k)\;,
$$ 
for all $k \geq 1$.~$\square$ 

\subsection{Proof of Corollary~\protect{\ref{relaxccc}}}
\label{proofcor}
\paragraph{}
To prove \reff{le.10}, first notice that by construction 
the process $(\widetilde U_n)_{n\in\rel}$ has the same law as
$(Z^{x}_n)_{n\in\rel}$ and $(\widetilde V_n)_{n\in\rel}$ has the
same law as $(Z^{y}_n)_{n\in\rel}$.  Thus,
\begin{equation}
\Bigl| \pro(Z^{x}_{n}= a) - \pro(Z^{y}_{n}= a) \Bigr| 
\;= \; \left| \pro(\widetilde{U}_n= a) - \pro(\widetilde{V}_n= a)
\right|
% \nonumber \\
%& = & \left| {\bf E}_{\pro}\left[1_{a}(\widetilde{U}_n) -
%  1_{a}(\widetilde{V}_n)\right] \right| \nonumber  \\
%& \leq & {\bf E}_{\pro}\left[1_{ \{ \widetilde{U}_n \neq
%\widetilde{V}_n \}} \ \right] \nonumber  \\ 
\; \leq \; \pro(\widetilde{U}_n \neq \widetilde{V}_n)) 
\label{premiere}
\end{equation}
Hence, by definition of the process $T_n$ and Lemma \ref{dominlem},
\begin{equation}
\Bigl| \pro(Z^{x}_{n}= a) - \pro(Z^{y}_{n}= a) \Bigr|
\;\leq\; %  \pro(\widetilde{U}_n \neq \widetilde{V}_n)) =  
\pro(T_{n} = 0) \leq  \pro(S^{(\gamma)}_n = 0)\;.
\label{le.30}
\end{equation}

\paragraph{}
The proof of \reff{le.7} starts similarly: 
\bqy \
\lefteqn{\hspace{-3cm}\Bigl|
  \pro\Bigl((Z^{x}_{n}, \ldots, Z^{x}_{n+k})= (a_0, \ldots, a_k)\Bigr)  -
  \pro\Bigl((Z^{y}_{n}, \ldots, Z^{y}_{n+k})= (a_0,
  \ldots, a_k) \Bigr) \Bigr| } \nonumber \\
%&= & \left| \pro\Bigl((\widetilde{U}_{n}, \ldots,
%  \widetilde{U}_{n+k})= (a_0, \ldots, a_k) \Bigr) -
%  \pro\Bigl((\widetilde{V}_{n}, \ldots, \widetilde{V}_{n+k})= (a_0,
%  \ldots, a_k) \Bigr)\right| \nonumber \\
%& = & \left| {\bf E}_{\pro}\left[1_{(a_0,
%  \ldots, a_k)}(\widetilde{U}_{n}, \ldots, \widetilde{U}_{n+k})  -1_{(a_0,
%  \ldots, a_k)}(\widetilde{V}_{n}, \ldots, \widetilde{V}_{n+k})\right] \right|
%\nonumber  \\
%& \leq & \pro(\Delta^{c}_{n,n+k})  \nonumber \\
%& \leq & \pro(T_{n+k} \leq k+1)  \nonumber \\ 
& \leq & \pro(S^{(\gamma)}_{n+k} \leq k+1) . \nonumber 
\eqy
%The last inequality is due to Lemma \ref{dominlem}.
To conclude, we notice that,  
\begin{equation}
\pro(S^{(\gamma)}_{n} \leq k) \;=\;  \sum_{j=0}^{k} \pro(S^{(\gamma)}_{n} = j) 
\;=\; \sum_{j=0}^{k}   \left(\prod_{m=1}^{j-1}(1- \gamma_m) \right)
\pro(S^{(\gamma)}_{n-j} = 0)   \;. \; \square  
\label{le.33}
\end{equation}

\section{Proof of Theorem~\ref{normalizedth}}% (normalized function)}
\label{normalized}

\paragraph{}
The proof of Theorem~\ref{normalizedth} is based on the inequality
\begin{equation}
\left| \int f \circ T^n g d\mu - \int f  d\mu \int g  d\mu \right|
%&= & \left| \int f(x)  \int \left( \esp[g((Z^{x}_{n+j})_{j\leq-1})] -
%\esp[g((Z^{y}_{n+j})_{j\leq-1})]  
%  \right) d\mu(y) d\mu(x) \right|   \nonumber \\
%&\leq  & \sup_{x,y}{\left|\esp[g((Z^{x}_{n+j})_{j\leq-1})] -
%\esp[g((Z^{y}_{n+j})_{j\leq-1})]\right|}  
%||f||_1   \nonumber \\
%&\leq  & \sup_{x,y}{\left|\esp[g((\tilde{U}_{n+j})_{j\leq-1})] -
%\esp[g((\tilde{V}_{n+j})_{j\leq-1})]\right|}  
%||f||_1   \nonumber \\
\;\leq  \; ||f||_1
    \,\sup_{x,y}{\esp\left[\left|g((\tilde{U}^{x,y}_{n+j})_{j\leq-1}) - 
    g((\tilde{V}^{x,y}_{n+j})_{j\leq-1})\right|\right]}  \;,
\label{eq:105.1}
\end{equation}
which follows from \reff{interpretation} and the fact that 
$\bigl((\widetilde{U}^{x, y}, \widetilde{V}^{x, y})\bigr)_{n\in\rel}$
is a coupling between the chains with pasts $x$ and $y$, respectively.
An upper bound to the right-hand side is provided by
Proposition~\ref{existcoupling}.  We see that the transition
probabilities \reff{eq:20} satisfy condition \reff{p.3}, since
\begin{equation}
  \label{eq:100.1}
  \frac{\trans(a \,|\, x )}{\trans(a \,|\, y )} 
\;=\; e^{\phi(ax) - \phi( ay)} 
\;\geq\;   e^{- \var_{m+1}(\phi)}
\end{equation}
whenever $x,y \in\underline{A}$ are such that
$x \stackrel{m}{=}y$ for some $m \in \nat$. 
We can therefore apply Proposition~\ref{existcoupling} with
\begin{equation}
  \label{eq:101}
  \gamma_m \;=\; 1 - e^{-\var_{m+1}(\phi)}\;,
\end{equation}
which tends monotonically to zero if $\sum_{m\geq 1} \var_m(\phi) < +
\infty$.

\paragraph{}
To prove \reff{eq:x} we use the process $(T^{x,y}_n)_{n\in {\bf Z}}$ to
obtain the upper bound 
\begin{eqnarray}
\esp\left[\left|g((\widetilde{U}^{x,y}_{n+j})_{j\leq-1}) -
  g((\widetilde{V}^{x,y}_{n+j})_{j\leq-1})\right|\right]
& =  & \esp\left[\sum_{k=0}^{+\infty}
    1_{ \{ T^{x,y}_{n}= k \} } \left|g((\tilde{U}_{n+j})_{j\leq-1}) -
      g((\tilde{V}_{n+j})_{j\leq-1})\right|\right]  \nonumber \\
&\leq  & \sum_{k=0}^{+\infty}\,  %\sup_{\scriptstyle z
%\stackrel{\scriptstyle k}{=} t}{|g(z)-g(t)|}\,
\var_k(g)\,\pro(T^{x,y}_{n}= k)    \nonumber \\
%&\leq  & \sum_{k=0}^{+\infty}\, \var_k(g)\,
%\pro(T^{x,y}_{n}= k)    \nonumber \\
&\leq  & ||g||_{\phi}\,\sum_{k=0}^{+\infty}\, \var_k(\phi)\,
\pro(T^{x,y}_{n}= k)\;.
\label{eq:110}
\end{eqnarray}
Now, in order to use the bound \reff{compteur} of
Proposition~\reff{existcoupling} we resort to the monotonicity of the
variations of $\phi$:
\begin{eqnarray}
\sum_{k=0}^{+\infty}\, \var_k(\phi)\, \pro(T^{x,y}_{n}= k)
&\leq  & \sum_{k=0}^{n-1}   \var_{k}(\phi)\, \pro(T^{x,y}_{n}= k) 
+ \var_{n}(\phi) \sum_{k=n}^{+ \infty} \pro(T^{x,y}_{n}= k) 
 \nonumber \\
&=  & \sum_{k=0}^{n-1}  \var_{k}(\phi) \pro(T^{x,y}_{n-k}= 0) +
\var_{n}(\phi) \sum_{k=n}^{+ \infty} \pro(T^{x,y}_{0}= k -n ) \nonumber\\
&\leq  & \sum_{k=0}^{n}   \var_{k}(\phi)\, \pro(S^{(\gamma)}_{n-k}= 0)\;,
\label{eq:115}
\end{eqnarray}
uniformly in $x,y$.  The bound \reff{eq:x} follows from
\reff{eq:105.1}, \reff{eq:110}, \reff{eq:115} and the fact that
\begin{equation}
  \label{eq:300}
 \sum_{j=0}^{+ \infty} \pro(T^{x,y}_{0}= j ) \;=\; 1 \;=\;
\pro(S^{(\gamma)}_{0}= 0) \;.
\end{equation}

\paragraph{}
To prove \reff{eq:xx} we use the strong Markov poroperty of the
process $(S^{(\gamma)}_n)_{n\in\nat}$ to obtain
\beq
\label{relationbasique}
\pro(S^{(\gamma)}_{n} = 0) = \sum_{k=1}^{n} \pro( \tau =  k ) \, \pro(
S^{(\gamma)}_{n-k} =0 )\;,
\eeq
where
\begin{equation}
\tau  \;=\;  \inf \{ n > 0 ; S^{(\gamma)}_n = 0 \} \;.
\label{eq:107}
\end{equation}

\paragraph{}
We now use (\ref{relationbasique}) to bound the last line in \reff{eq:115}
in the form
\bqy
\sum_{k=0}^{n}   \var_{k}(\phi)\, \pro(S^{(\gamma)}_{n-k}= 0)&\leq  &
\sum_{k=1}^{n}\bigl[\var_0(\phi) \, \pro( \tau =  k ) +
\var_k(\phi)\bigr]\,  \pro(S^{(\gamma)}_{n-k}= 0) \nonumber \\
&\leq  &  C \sum_{k=1}^{n}   \pro(\tau = k)\,\pro(S^{(\gamma)}_{n-k}=
0) \nonumber \\
&= &  C \;  \pro(S^{(\gamma)}_{n}= 0) \;,
\label{eq:122}
\eqy
with
\begin{equation}
C \;=\; \var_0(\phi) + \sup_k \,{\var_k(\phi) \over \pro( \tau =
k )}\;.
\label{eq:122.1}
\end{equation}
To conclude, we must prove that the constant $C$ is finite.  By direct
computation,
\begin{eqnarray}
\label{tau}
\pro(\tau  = 1)  &= & \gamma_{0},\nonumber \\[3pt]    
\pro(\tau = n)& =& \gamma_{n-1}\,\prod_{m=0}^{n-2}(1- \gamma_{m})
\quad  \hbox{ for } n\geq 2,\\[3pt]
\pro(\tau = +\infty) &= & \prod_{m=0}^{+ \infty}(1-
\gamma_{m})\;. \nonumber
\end{eqnarray}
From this and \reff{rr.30} we obtain
\begin{equation}
\lim_{k\to\infty} {\var_k(\phi) \over \pro( \tau = k )} \;=\;
\lim_{k\to\infty} \,{\var_k(\phi) \over 1 - e^{-\var_k(\phi)}}\,
{1 \over \prod_{m=0}^{k-2}(1- \gamma_{m})} \;.
%\;=\; {1 \over \pro(\tau = +\infty)}\;,
\end{equation}
Since $\var_k(\phi)\to 0$, the first fraction converges to 1.  We see
from \reff{tau} that the second fraction converges to 
$1/\pro(\tau = +\infty)$.  By elementary calculus, this is finite
since $\phi$ has summable variations.~\square

\begin{rem}
The previous computations lead to stronger results for more regular
functions $g$. For example, when $g$ satisfies
\begin{equation}
  \label{eq:301}
  \var_k(g) \;\le\; ||g||_{\theta} \,\theta^k
\end{equation}
for some $\theta <1$ and some $||g||_{\theta}<\infty$ (H\"older norm
of $g$), a chain of inequalities almost identical to those ending in
\reff{eq:110} leads to 
\bqy \left| \int f \circ T^n g d\mu - \int f
  d\mu \int g d\mu \right| &\leq& ||f||_1\,\sum_{k=0}^{+\infty}
||g||_{\theta}\, \theta^{k}
\,\gamma^*_{n-k} \nonumber \\
&\leq & ||f||_1 \,||g||_{\theta}\,\theta^n \sum_{k=0}^{n}
\theta^{-k}\,\gamma^*_{k} \;.
\label{eq:122.2}
\eqy 

On the other hand, if $g$ is a function that depends only on the first
coordinate, we get, 
\bqy 
\left| \int f \circ T^n g d\mu - \int f d\mu \int g d\mu
\right| &\leq & ||f||_1 \,\sup_{x,y}{\Bigl|\esp[g(Z^{x}_n)] -
\esp[g(Z^{y}_n)]\Bigr|}\nonumber \\ 
&\leq & ||f||_1\, ||g||_{\infty} \pro( \tilde{U}_n \neq \tilde{V}_n)
\nonumber \\ 
&\leq & ||f||_1\, ||g||_{\infty} \, \gamma^*_n\;.  
\label{eq:130}
\eqy
\end{rem}

\section{Proof of Theorem~\ref{general}}
\label{proof}
\paragraph{}
We now consider the general case where the function $\phi$ is not
necessarily normalized.  In this case we resort to the normalization 
$\psi$ define in \reff{rr.25} and we consider
chains with transition probabilities
\begin{equation}
\trans(a \,|\, x ) \;=\; e^{\phi(xa)}\,\frac{\rho(xa)}{\rho(x)}
\;=:\; e^{\psi(xa)} \;.
\label{gra.10}
\end{equation}
However, the summability of the variations of $\phi$ does not imply
the analogous condition for $\psi$, because there
are addition ``oscillations'' due to the cocycle $\log{\rho} -
\log{\rho \circ T}$.  Instead,  
\begin{equation}
\left.
\begin{array}{r}
\var_m\psi\\
\var_m(\log\rho)
\end{array}\right\} \;\leq\; \sum_{k \geq m} \var_k(\phi)\;,
\label{gra.0}
\end{equation}
for all $m\ge 0$ (see Walters 1978).
Hence, we can apply Theorem~\ref{normalizedth} only under the
condition 
\begin{equation}
\sum_{k=1}^{+\infty} k\, \var_k(\phi) < + \infty\;.
\label{gra.1}
\end{equation}
If this is the case, the correlations for functions $f \in L^1(\mu)$
and $g \in V_\psi$ decay faster than $\gamma^*_{m}$, where $\gamma_{m}
= e^{\sum_{k \geq m} \var_k(\phi)} -1$.
%If $\var_m(\phi)$ decreases (at most) polynomially, then the speed
%of the decrease of correlations is of order $\sum_{k\geq m}\var_k(\phi)$. 

\paragraph{}
To prove the general result without assuming \reff{gra.1} we must work
with \emph{block} transition probabilities, which are less sensitive
to the oscillations of the cocycle.  More precisely, 
given a family of transition probabilities $\trans$ on $A \times
\underline{A}$, let $\trans_n$ denote the corresponding transition
probabilities on $A^n \times \underline{A}$:
\begin{equation}
\trans_{n+1}(\atom{a}{0}{n} \,|\, x) %= \trans(a_n | a_{n-1} \cdots
%a_1 x) \trans_{n}(\atom{a}{0}{n-1}| x) = 
\;=\; \trans(a_n \,|\, a_{n-1} \cdots\, a_1 x) \cdots 
\trans(a_2 \,|\,a_1 x)\, \trans(a_1 \,|\, x)
\label{gra.3}
\end{equation}
where 
\begin{equation}
\atom{a}{0}{n} \;:=\; (a_0, \ldots, a_n) \;\in A^{n+1}\;. 
\label{gra.4}
\end{equation}
If the transition probabilities $\trans$ are defined by a normalized
function $\phi$ as in \reff{eq:20}, then we see from \reff{gra.3} that
the transition probabilities $\trans_n$ obey a similar relation
\begin{equation}
  \label{gra.7}
  \trans_n(\atom{a}{0}{n-1}\, |\, x ) = e^{\phi_n(x\atom{a}{0}{n-1})}\;,
\end{equation}
with
\begin{equation}
\phi_n(x\atom{a}{0}{n-1}) \;:=\; \sum_{k=0}^{n-1}\phi(x a_0 \cdots
a_k)\;.
\label{gra.8}
\end{equation}
In particular, for transitions \reff{gra.10} the formula \reff{gra.3}
yields
\begin{equation}
\psi_n \;=\; \phi_n + \log{\rho} - \log{\rho \circ T^n} \;.
\label{gra.12}
\end{equation}

\paragraph{}
A comparison of \reff{gra.12} with \reff{gra.0} shows that it is
largely advantageous to bound directly the oscillations of $\psi_n$.
This is what we do in this section by adapting the arguments of
Section~\ref{normalized}.

\subsection{Coupling of the transition probabilities for blocks}
\paragraph{}
For every integer $n$, we define a family of transition probability 
$\transb_n$ on $(A^{n})^2 \times \underline{A}^2$ by 
\begin{equation}
\transt_n(\atom{a}{0}{n-1}, \atom{b}{0}{n-1}\,|\, x, y) \;=\;
\bigl[ \trans_n(\cdot \,|\,x)\, \tilde{\times}\, \trans_n(.  \,|\,
y)\bigr]
(\atom{a}{0}{n-1}; \atom{b}{0}{n-1})\;.
\label{gra.27}
\end{equation}

Let $(n_{m})_{m \in \nat}$ be an increasing sequence. 
For each double history $x,y$, we consider the coupling 
$\bigl((\overline{U}^{x, y}, \overline{V}^{x, y})\bigr)_{m\in\rel}$
of the chains for $n_m$-blocks with past $x$ and $y$,  defined by, 
\begin{eqnarray}
\lefteqn{
\pro(\atom{\overline{U}^{x,y}}{0}{n_m}= \atom{a}{0}{n_m}\,,\,
\atom{\overline{V}^{x,y}}{0}{n_m}= \atom{b}{0}{n_m})}
\nonumber\\
&&=\ \prod_{m=1}^{M} 
\transb_{n_{m+1} - n_{m}}(\atom{a}{n_m}{n_{m+1}} \,,\, \atom{b}{n_m}{n_{m+1}}
\,|\, a_{n_m} \cdots a_0 x\,,\, b_{n_m} \cdots b_0 y)\;.
\label{gra.31}
\end{eqnarray}

\subsection{The process of last block-differences}
\paragraph{}
We set 
\begin{equation}
\gamma^{(n)}_k \;=\; 1 - \inf\, \left \{
{\frac{\trans_n(\atom{a}{0}{n-1}\,|\, x) }{\trans_n(\atom{a}{0}{n-1}
\,|\, y)}}
\,:\, x \stackrel{k}{=}y\,,\; a_1, \ldots, a_{n-1} \in A\right\}\;.
\end{equation}
From~(\ref{coupling2}) we see that, for $x \stackrel{k}{=}y$, the
weight of the diagonal of each coupling $\transb_n$ satisfies
\beq
\label{diago}
\Delta\bigl(\transt_n(\cdot,\cdot \, |\, x, y)\bigr) \; \geq\;
 \inf_{a_0, \ldots, a_{n-1} \in
 A}{\frac{\trans(\atom{a}{0}{n-1}\,|\, x)
 }{\trans(\atom{a}{0}{n-1} \,|\, y)}}\; \geq \; 1 - \gamma^{(n)}_k \;.
\eeq
If we denote %by $\overline{\Delta}^{x,y}_{m, m+q}$ the set 
$$\overline{\Delta}^{x,y}_{m, m+q} \;:= \;
%\bigcap_{n_m \leq j \leq n_{m+q} } %\{ \widetilde{U}_j =
%\widetilde{V}_j \} =  
\Bigl\{ \overline{U}^{x,y}_j =  \overline{V}^{x,y}_j \, , \,  n_m \leq
j \leq n_{m+q}\Bigr\} \;,$$
we deduce from \reff{diago} that
\beq
\label{onestep}
\pro(\overline{\Delta}_{m+k+1} \,|\, \overline{\Delta}_{m,m+k})
\;\geq\; 1 - \gamma^{(n_{m+k+1}-n_{m+k})}_{n_{m+k} - n_m} \;.
\eeq

\paragraph{}
We construct the process $(\overline{T}_n)_{n \in \nat}$ with
\begin{equation}
\overline{T}^{x,y}_m \;=\; \inf \,\Bigl\{ p \geq 0 \,:\, 
U^{x,y}_i \neq V^{x,y}_i \hbox{ for some } i, n_{m-p} \leq i \leq
n_{m-p+1} \Bigr \} \;.
%&=& \inf \,\Bigl \{ p \geq 0 \,:\, \widetilde{Z}^{x,y} \not \in
%\overline{\Delta}_{m-p, m-p+1}\} $$ 
\end{equation}
By \reff{onestep}, the conditional laws of this process satisfy, 
\begin{equation}
\pro(\overline{T}_{m+1}=k+1 \,|\, \overline{T}_{m} = k) \;\geq\;
  1 - \gamma^{(n_{m+k+1}-n_{m+k})}_{n_{m+k} - n_m} 
\label{gra.40}
\end{equation}
and 
\begin{equation}
\pro(\overline{T}_{m+1}= 0 \;|\; \overline{T}_{m} = k) \;\leq\;
\gamma^{(n_{m+k+1}-n_{m+k})}_{n_{m+k} - n_m}\;. 
\label{gra.41}
\end{equation}

\subsection{The dominating Markov process}
\paragraph{}
Let us choose the length of the blocks in such a way that the sequence
$(n_{m})_{m \in \nat}$ is subadditive, i.e.
\begin{equation}
n_{m+k} - n_m \;\leq\; n_k
\label{gra.sub}
\end{equation}
for $m,k\ge 0$, and that
\begin{equation}
\sup_{n\ge 0} \gamma^{(n)}_{\ell} \;<\; 1
\label{gra.sup}
\end{equation}
for all $\ell\ge 0$.
These two properties together with \reff{gra.40}--\reff{gra.41} imply
that, for all histories $x$ and $y$,
\begin{equation}
\pro(\overline{T}^{x,y}_{m+1}=k+1 \,|\, \overline{T}^{x,y}_{m} = k) \;\geq\;
  1 - \overline{\gamma}_{k} 
\label{gra.50}
\end{equation}
and 
\begin{equation}
\pro(\overline{T}^{x,y}_{m+1}= 0 \,|\, \overline{T}^{x,y}_{m} = k) \;\leq\;
\overline{\gamma}_{k}. 
\label{gra.51}
\end{equation}
with
\begin{equation}
\overline{\gamma}_{k} \;:=\; \sup_{n \geq 1} \gamma^{(n)}_{n_k}\;,
\label{gra.52}
\end{equation}
for $m \geq 1$.

\paragraph{}
We now define the ``dominating'' Markov chain $(S^{(\overline\gamma)}_n)_{n \in
  \nat}$ as in \reff{rr.5.-1}--\reff{rr.5}. Lemma~\ref{dominlem} yields 
\begin{equation}
\pro(\overline{T}^{x,y}_m = 0) \;\leq\; \pro(\overline{S}_m = 0)
\;\leq\; \overline{\gamma}^*_m\;.
\end{equation}
Hence, if $n_m \leq n \leq n_{m+1}$,
\begin{equation}
\pro(\overline{U}^{x,y}_n \neq \overline{V}^{x,y}_n) \;\leq\;
\pro(T^{x,y}_m = 0) \;\leq\;  \overline{\gamma}^*_m\;.
\end{equation}

\subsection{Decay of correlations}

\paragraph{}
We can now mimick the proof of Theorem~\ref{normalized} in terms of
barred objects.  

\paragraph{}
As $( \var_m(\phi))_{m \in \nat}$ is summable, there exists a
subadditive sequence $(n_{m})_{m \in \nat}$ such that the sequence
$\alpha_m$ of the tails
%$(\sum_{k \geq n_{m}}  \var_k(\phi))_{m \in \nat}$ 
\begin{equation}
\alpha_m \;=\; \sum_{k \geq n_{m}}  \var_k(\phi)
\end{equation}
is summable:
\begin{equation}
\sum_{m \geq 0}  \alpha_m < + \infty\;.
\label{gra.51.1}
\end{equation}

\paragraph{}
The transitions for blocks of size $n$ satisfy
\begin{equation}
\frac{\trans_n(\atom{a}{0}{n-1} \,|\, x) }{\trans_n(\atom{a}{0}{n-1}
  \,|\, y)} \;\ge\; e^{-\var_k(\psi_n)}
\label{gra.55}
\end{equation}
if $x\stackrel{k}{=} y$. 
But from \reff{gra.12}, \reff{gra.8} and \reff{gra.0} we have
\begin{eqnarray}
\var_k(\psi_n) &\le& \left(\sum_{m=k}^{k+n}  + \sum_{m\ge k+n} +
\sum_{m\ge k}\right) \var_m(\phi) \nonumber\\
&\le& 3\, \sum_{m\ge k}\var_m(\phi)\;.
\label{gra.56}
\end{eqnarray}
Hence we can choose in \reff{gra.52}
\begin{equation}
\overline{\gamma}_{k}  \;\leq\; 1 - e^{- 3 \alpha_k}\;,
\label{gra.57}
\end{equation}
a choice for which
\begin{equation}
\sum_{k\geq 1} \overline{\gamma}_{k} < + \infty \;.
\label{gra.58}
\end{equation}

\paragraph{}
To prove the theorem, we now proceed as in \reff{eq:105.1} and
\reff{eq:110}--\reff{eq:122.1} but replacing tildes by bars and
putting bars over the processes $(T_n)$ and $(S_n^{(\gamma)})$.  We just
point out that, due to the subadditivity of $n_m$,
$$\var_{(n_{m+k} - n_m)}(\phi) \; \leq\; \var_{n_k}(\phi)$$
uniformly in $m$.~\square

%\begin{rem}
%If  $\overline{\gamma}^*_m = O(\alpha_m)$, noticing that $\alpha_m
%= \sum_{k \geq n_m
% }  \var_k(\phi)$, we get 
%\footnote{One possibility to get sure to be in this case is to choose the subsequence in
% such way that $\alpha$ has (at most) polynomial decreasingy. This should be possible because
% we are in the case when $\gamma_m$ decreases slowly (rest not summable). We
% must investigate this carefully.} 
%$$\pro(\tilde{U}_n \neq \tilde{V}_n) = O\left(\sum_{k \geq n
% } \var_k(\phi)\right).  $$
%\end{rem}

\appendix

\section{Returns to the origin of the dominating Markov chain}
\label{markovchain}
\paragraph{}
In this appendix we collect a few results concerning the
probability of return to the origin of the Markov chain
$(S^{(\gamma)}_{n})_{n \in \nat}$ defined via \reff{rr.5}.  
(In the sequel we omit the superscript ``$(\gamma)$'' for simplicity.)
 
\begin{prop}
\label{asymptotic}
Let $(\gamma_{n})_{n \in \nat}$ be a real-valued sequence decreasing
to $0$ as $n \to + \infty$. 
\bid
\item[(i)] If $ \displaystyle \sum_{m\geq 1} \prod_{k=0}^{m}(1 -
  \gamma_k) = + \infty $, then $\pro(S_n = 0) \to 0$.
\item[(ii)] If $ \displaystyle \sum_{m\geq 1} \gamma_k < + \infty $,
  then $\sum_{n \geq 0} \pro(S_n = 0) < + \infty$.
\item[(iii)] If $(\gamma_m)$ decreases exponentially, then so does
  $\pro(S_n = 0)$.
\item[(iv)] If $(\gamma_m)$ decreases polynomially, then $\pro(S_n = 0)
  = O(\gamma_n)$.
\eid
\end{prop}

\paragraph{Sketch of the proof}

\paragraph{}
Statement $(i)$ follows from the well known fact that the
Markov chain $(S_n)_{n\in\nat}$ is positive recurent if and only if, 
$$\sum_{m\geq 1} \prod_{k=0}^{m}(1 - \gamma_k) <  + \infty. $$

\paragraph{}
To prove parts (ii) and (iii) we introduce the series
\begin{equation}
F(s) = \sum_{n=1}^{+ \infty} \pro(\tau=n) \,s^{n}\;,
\label{gra.61}
\end{equation}
and
\begin{equation}
G(s) \;=\; \sum_{n=0}^{+ \infty} \pro(S_n =0)\, s^{n} 
\label{gra.60}
\end{equation}
where the random variable $\tau$ is the time of first return to zero,
defined in \reff{eq:107}.  The probabilities $\pro(\tau=n)$ were
computed in \reff{tau} above.
The relation (\ref{relationbasique})
implies that these series are related in the form
%\beq
%\label{avecH}
%G(s) H(s) = \frac{1}{1 -s }
%\eeq
%for all $0 \leq s<1$, or, after a simple computation, 
\beq
\label{avecF}
G(s) = \frac{1}{1 -F(s)}\;,
\eeq
for all $s \geq 0$ such that $F(s)<1$.  

\paragraph{}
It is clear that the radius of convergence of $F$ is at least 1.  In fact, 
  \begin{equation}
    \label{eq:151}
    F(1) \;=\; \pro(\tau < +\infty)\;.
  \end{equation}
Moreover, if $ \sum_{m\geq 1} \gamma_k < + \infty $, the radius of
  convergence of $F$ is
  \begin{equation}
    \label{eq:152}
    \lim_{n\to\infty} \, [\gamma_n]^{-1/n}\;.
  \end{equation}
This is a consequence of the fact that
$\pro(\tau=n)/\gamma_{n-1}\to \pro(\tau=+\infty) > 0$, as concluded
from \reff{tau}.

\paragraph{}
Statement $(ii)$ of the proposition is a consequence of the fact that
the radius of convergence of the series $G$ is at 
  least 1 if $ \sum_{m\geq 1} \gamma_k < + \infty $.  This follows
  from the relation \reff{avecF} and the fact that the right-hand side
  of \reff{eq:151} is strictly less than one when the chain
  $(S_n^{(\gamma)})$ is transient.

\paragraph{}
To prove statement (iii) let us assume that $\gamma_m \leq C \gamma^m$
for some constants $C < +\infty $ and $0<\gamma<1$.  By \reff{eq:152},
the radius of convergence of $F$ is $\gamma^{-1} > 1$ while, by
\reff{eq:151}, $F(1)<1$.  By continuity it follows that there exists
$s_0>1$ such that $F(s_0)=1$ and, hence, by \reff{avecF},
$G(s)<+\infty$ for all $s<s_0$.  By definition of $G$, this implies
that $\pro(S_n = 0)$ decreases faster than $\zeta^n$ for any
$\zeta\in(1,s_0^{-1})$.  

\paragraph{}
Statement (iv) is a consequence of the following lemma.

\begin{lem}
If 
\beq
\label{condpoly}
\alpha \;:=\; \sup_i\, \overline{\lim}_{k\to\infty}
\left[ \frac{\pro(\tau = i)}{\pro(\tau =ki)}\right]^{1/k}
\;<\; {1\over \pro(\tau<+\infty)}\;,
\eeq
then
$$  \pro(S_n = 0) = O\left( \pro( \tau = n) \right).$$
\end{lem}

\paragraph{Proof}
We start with the following observation.  If 
$i_1 + \cdots + i_k = n$, then $\max_{1 \leq m \leq k} i_m . n/k$ and
thus, for $g$ is an increasing
$$ 
g(n) \;\leq\; g\left( k\, i_{\rm max}\right) \;,
$$
where $i_{\rm max} = \max_{1\le m\le k}\,i_m$.  
If we apply this to $g(n) = 1/\pro(\tau=n)$, which is increasing by
\reff{tau}, we obtain
\begin{equation}
1 \;\le\; \frac{\pro( \tau = n)}{\pro(\tau = k\,i_{\rm max})}\;.
\label{fra.1}
\end{equation}

\paragraph{}
We now invoke the following a explicit relation between the
coefficients of $F$ and $G$.
\beq
\label{generalexpression}
\pro(S_n =0) = \sum_{k = 1}^{n} \sum_{ \scriptsize 
\begin{array}{c}
i_1, \ldots, i_k \geq 1 \\
i_1 + \cdots + i_k  = n 
\end{array}
} \prod_{m=1}^{k} \pro(\tau = i_m) \;,
\eeq
for $n \geq 1$.  
%This can be obtained  directly
%from~(\ref{avecF}) or, alternatively, by
%decomposing each return time as a sum of times of \emph{first}
%returns.  
Multiplying and dividing each factor in the rightmost
product by $\pro(\tau<+\infty)$, this formula can be rewritten as 
\beq
\label{autreexpression}
\pro(S_n =0) = \sum_{k = 1}^{n} \pro(\tau < + \infty)^k  \sum_{ \scriptsize 
\begin{array}{c}
i_1, \ldots, i_k \geq 1 \\
i_1 + \cdots + i_k  = n 
\end{array}
} \prod_{m=1}^{k} \pro(\tau = i_m \,|\, \tau <
      +\infty). 
\eeq 
Combining this with \reff{fra.1} we obtain
\begin{equation}
\label{fra.2}
\pro(S_n =0) \;\le\; \pro(\tau=n)\,\sum_{k = 1}^{n} \pro(\tau < +
\infty)^k  \sum_{ \scriptsize  
\begin{array}{c}
i_1, \ldots, i_k \geq 1 \\
i_1 + \cdots + i_k  = n 
\end{array}
} \prod_{m=1}^{k} {\pro(\tau = i_m \,|\, \tau < +\infty)\over 
\pro(\tau = k\, i_{\rm max})}\;,
\end{equation}
If we single out the factor $\pro(\tau=i_{\rm max}\,|\, \tau <
+\infty)=\pro(\tau=i_{\rm max})/\pro(\tau<+\infty)$ from the rightmost product
of \reff{fra.2} and use the hypothesis \reff{condpoly} we get
\begin{equation}
\label{fra.3}
\pro(S_n =0) \;\le\; C\,\pro(\tau=n)\,\sum_{k = 1}^{n} \alpha^k\,\pro(\tau < +
\infty)^{k-1}  \sum_{ \scriptsize  
\begin{array}{c}
i_1, \ldots, i_k \geq 1 \\
i_1 + \cdots + i_k  = n 
\end{array}
} \prod_{\scriptsize 1\le m \le k\atop \scriptstyle i_m\neq i_{\rm max}} 
\pro(\tau = i_m \,|\, \tau < +\infty)\;,
\end{equation}
for some constant $C>0$.
To bound the last sum on the right-hand side 
we introduce a sequence of independent random variables
$(\tau^{(i)})_{i \in \nat}$ with common distribution
\begin{equation}
\pro(\tau^{(i)} = j) \;=\; \pro(\tau = j \,|\, \tau < +\infty)\;.
\label{fra.0}
\end{equation}
Then
\begin{equation}
\sum_{j= 1}^{n-k+1} \pro\Bigl( \sum_{i=1}^{k-1} \tau^{(i)} =
n-j\Bigr) \;\le\; 1\;.
\label{fra.5}
\end{equation}
Hence, \reff{fra.3} implies
\begin{equation}
\label{fra.4}
\pro(S_n = 0) \;\leq\; C\, \alpha\, \sum_{k=1}^{\infty} 
\left[\alpha^{-1}\pro(\tau<+ \infty)\right]^{k-1}\,
\pro(\tau = n) \;\leq\; \hbox{const}\, \pro(\tau = n)\;. \; \square 
\end{equation}

\paragraph{}
We notice that, according to~(\ref{tau}), $\gamma_n \sim 
\pro(\tau = n)/\pro(\tau=+ \infty)$. Hence, a sufficient condition for
\reff{condpoly} is a similar condition for the sequence
$(\gamma_n)$.  Such a condition holds, for instance, if the later
sequence decays polynomially. Statement $(iii)$ of the 
proposition follows.~\square

\subsection*{Acknowledgements}  It is a pleasure to thank Pablo Ferrari
and Luis Renato Fontes for useful discussions.  We also thank Pablo
Ferrari for informing us about his work in progress with Alejandro
Maass and Servet Mar\'{\i}nez on a related regenerative representation
of chains with complete connections.

\vskip30pt

Xavier Bressaud and Antoni Galves

Instituto de Matem\'atica e Estat\'{\i}stica

Universidade de S\~ao Paulo

Caixa Postal 66281

05315-970 S\~ao Paulo, Brasil

e-mail: {\tt bressaud, galves@ime.usp.br}

\bigskip

Roberto Fern\'andez

Instituto de Estudos Avan\c{c}ados

Universidade de S\~ao Paulo

Av.\ Prof.\ Luciano Gualberto, Travessa J, 374 T\'erreo,

Cidade Universit\'aria

05508-900 S\~ao Paulo, Brasil

e-mail: {\tt rf@ime.usp.br}

\end{document}